\documentclass{amsart}

 \def\newblock{\ }%
\usepackage{amsmath,verbatim,mathtools}
\usepackage{graphicx}
\usepackage{color}
\usepackage{psfrag}
\usepackage[makeroom]{cancel}
\newcommand{\al}{\alpha}                
\newcommand{\g}{\lambda}                
\newcommand{\Cplus}{{\mathfrak C}_+(y;\g)}
\newcommand{\Cminus}{{\mathfrak C}_-(y;\g)}
\newcommand{\Bplus}{{B}_+(y)}
\newcommand{\BMminus}{{B}^{(m)}_-(y)}
\newcommand{\Bdminus}{{B}^{(d)}_-(y)}
  
\newcommand{\R}{\mathbb{R}}
\newcommand{\Z}{\mathbb{Z}_{\ge 0}}
\DeclareTextFontCommand{\textbfit}{%
  \fontseries\bfdefault 
  \itshape
}

\newtheorem{proposition}{Proposition}
\newtheorem{theorem}{Theorem}[section]
\newtheorem{lemma}{Lemma}[section]
\newtheorem{remark}{Remark}[section]
\newtheorem{asm}{Assumption}
\newtheorem{cor}{Corollary}

\newtheorem{example}{Example}

\begin{document}




\title[Relative Error of Scaled Poisson]{Relative Error of Scaled Poisson Approximation via Stein's Method}


\author{ Yingdong Lu}
\address{IBM Research}
\email{yingdong@us.ibm.com}

\author{Yue Tan}
\address{Ohio State University}
\email{tan.268@buckeyemail.osu.edu}

\author{Cathy H. Xia}
\address{Ohio State University}
\email{xia.52@osu.edu}



\maketitle

\begin{abstract}
We study the accuracy of \emph{a scaled Poisson approximation} to the weighted sum of independent Poisson random variables, focusing on in particular the relative error of the tail distribution.
A bound on the relative approximation error is established using a modified Stein-Chen method. 
\end{abstract}%

%

\section{Introduction} 
\label{sec:intro}

Weighted summation of independent Poisson random variables (r.v.s) is a probabilistic entity that plays a crucial role in a wide variety of applications, such as epidemiology~\cite{DKES}, physics~\cite{BohmZech}, and computer science~\cite{MAMA2012}. 
A traditional approach of 
its quantification 
is a normal approximation matching the first two moments, which can be effective when the mean is large. In \cite{Fay}, a \emph{scaled} Poisson approximation that matches the first two moments after a scaling is proposed. 
That is, the distribution of $S$, 
a weighted sum of independent Poisson random variables
with mean $\mu = E[S]$ and variance $\sigma^2 = Var(S)$,
is approximated by that of a scaled Poisson random variable, i.e.,
\begin{equation} \label{eqn:PoisApprox}
S\approx_d \frac{1}{\beta}A_{\beta\mu}, \quad \quad 
\quad \beta = \frac{\mu}{\sigma^2},
\end{equation}
where $A_{\lambda}$ denotes a Poisson r.v. with mean $\lambda$. 
This way, the first two moments of $S$ are matched by that of the approximation. 
Through numerical experiments,  \cite{Fay} reports that this approach produces a more accurate tail approximation than the normal approximation. Moreover, the approximation is often more \emph{conservative} which serves certain types of application well, see, e.g. \cite{Fay,Ng2008,
Sig2012}.
However, no theoretical analysis on the quality of this approximation has been given nor is the approximation error well understood. It is the purpose of this paper to address these issues 
by providing quantitative bounds on the relative error of this approximation,
utilizing a novel application of the Stein method. 

The Stein's method has been a powerful tool for bounding the error of an approximating distribution to an unknown distribution.
The method was first introduced by Charles Stein~\cite{Stein1972} to study the normal approximation to the distribution of a sum of dependent r.v.'s.
It was later extended to handle Poisson approximations in~\cite{Chen1975}, often referred to as the Stein-Chen method specifically for studying Poisson approximations.
The Stein's method has been used to study many approximating distributions including binomial distribution \cite{Ehm1991}, multinomial distribution \cite{Loh1992},
compound Poisson distribution \cite{Barbour1992}, Gamma distribution \cite{Luk1994}, geometric distribution \cite{Pekoz1996}, etc. 
Applications in other contexts include ~\cite{chatterjee2010}, \cite{Ying:2016:AEM:2964791.2901463}, 
~\cite{Liu:2016:KSD:3045390.3045421},~\cite{Chwialkowski:2016:KTG:3045390.3045665},
~\cite{braverman2016}, and~\cite{braverman2017}.

The Stein's method relies on identifying a Stein's operator $\mathcal{A}$ associated with a given distribution~$Q$ such that a r.v. $Y$ follows distribution $Q$ (written as $Y \sim Q$) if and only if 
$E(\mathcal{A}f)(Y) =0$ for all bounded real-valued functions $f$ defined on the range of $Y$. If the distribution of a r.v.~$W$ can be well approximated by $Q$, then
$E(\mathcal{A}f)(W)\approx 0$.  
The error of approximating $E[h(W)]$ by $E[h(Y)]$, for a given metric function $h$, 
can then be estimated by studying the behavior of $E\mathcal{A} f_h(W)$.
where $f_h$ is a solution to the so-called Stein's equation:
\begin{equation}\label{eqn:SteinNorm}
\mathcal{A}f(w) = h(w)-Eh(Y).
\end{equation}


The
Stein's operator $\mathcal{A}$ frequently used in studying Poisson approximations was provided by
~\cite{Chen1975} and often referred to as the Stein-Chen method:
\begin{equation} \label{eqn:SteinOp1}
\mathcal{A}f(w) = \lambda 
\cdot  f(w+1) -w f(w), 
\end{equation}
%
Using the Stein-Chen method, moderate deviation bounds on the relative error of a Poisson approximation 
to a summation of locally dependent indicators 
was established in~\cite{ChenFangShao2013},
and were later made more general and explicit in~\cite{liu_xia_2020}. 

For the scaled Poisson approximation~\eqref{eqn:PoisApprox}, the scalar $\beta=\frac{\mu}{\sigma^2}$ is typically not integer-valued, 
thus the conventional Stein's operator $\mathcal{A}$ given by~\eqref{eqn:SteinOp1}
\textcolor{black}{can not be applied directly} as it only acts on bounded functions $f$ defined on the integers.
We show that, when the scalar is a rational number with $\beta=\frac{n}{m}$,
the scaled Poisson approximation  \eqref{eqn:PoisApprox}, now equivalent to $nS \approx m A_\lambda$ with $\lambda = \beta \mu$, 
can still be analyzed via the Stein's method. 
We introduce
\textcolor{black}{a new Stein's operator for 
$\hat{A}_\lambda:= m A_\lambda$,
defined as:
\begin{equation} \label{eqn:SteinOp2}
\mathcal{A}f(w) =  \lambda m f(w+m)-wf(w),
\end{equation} 
which provides a more suitable framework for  studying the approximation.}
This requires us to obtain new solutions to this Stein's equation 
and study the quantitative properties of the solution,
which form the most technical part of the paper. These properties 
allow us to apply  Stein's method and derive bound on the relative error of scaled Poisson approximation. 
\textcolor{black}{
We then extend our relative error bound to the case when $\beta$ is irrational using bounded convergence theorem.}

We emphasize that while the scaled Poisson approximation has been widely used in practice, to the best our knowledge,  this is the first work to quantify the relative error of the approximation via moderate deviations bound.
Meanwhile, the scaled Stein method provides a new method for analyzing approximations on a grid with applications beyond the weighted Poisson summation discussed in this paper. Moreover, the techniques developed in estimating differences of the solutions to the Stein's equations appeared to reveal some fundamental properties of Poisson distribution function, should be of independent interests.

Our bound on the relative error is similar in spirit and numerically comparable to those obtained in~\cite{liu_xia_2020} for (non-scaled) Poisson approximation.
Furthermore, we demonstrate that one of the conjectures raised in~\cite{liu_xia_2020} is valid beyond a threshold for which we provide an explicit expression. 

The rest of the paper is organized as follow. In Section \ref{sec:results}, we present the main results of the paper, which include a new Stein's operator, 
our 
relative error bound for the scaled Poisson approximation, 
and a proposition to address one of the conjectures raised in~\cite{liu_xia_2020}. 
Proof for the main results is presented in Section \ref{sec:proofs}; Key properties of the solution to modified Stein's equation are collected in Section~\ref{sec:properties}. 
The paper is concluded in Section~\ref{sec:conclusions} with summary and ongoing research. 

\section{Main Results and Implications}
\label{sec:results}

In this section, we present the problem of scaled Poisson approximation, the main results and their implications. Assumptions of the scaled Poisson approximation are presented in Sec. \ref{sub:1};  a modified Stein's method is introduced in Sec. \ref{sec:modifiedStein}; and our moderate deviations result is elaborated in Sec. \ref{sec:moderateResults}.


\smallskip
\subsection{Scaled Poisson Approximation}\label{sub:1}

Consider a weighted summation of $R$ independent Poisson r.v.s $\{A_{\nu_r}\}_{r=1}^R$,
\begin{equation}
S:=\sum_{r=1}^R b_rA_{\nu_r}, \label{eqn:P1}
\end{equation}
where  $A_{\nu}$ denotes a Poisson random variable with mean $\nu$ 
and $\{b_r\}_{r=1}^R$ are a set of positive weights. 
Without loss of generality, we assume that $b_r$'s are distinct and $0< b_1<b_2<\cdots<b_R$. If there exist $r$ and $s$ such that $b_r=b_s$, 
one can simply merge the two classes into one as $b_rA_{\nu_r}+b_sA_{\nu_s}\overset{d}{=}b_rA_{\nu_r+\nu_s}$.

Denote the mean and variance of $S$ as $\mu = E[S]$, and $\sigma^2 = Var(S)$ respectively, thus, 
\begin{equation} \label{eqn:parameters}
\mu = \sum_{r=1}^Rb_r\nu_r, \quad \quad \sigma^2=\sum_{r=1}^Rb_r^2\nu_r.
\end{equation} 
Let $\beta = \frac{\mu}{\sigma^2}$, which is known as the inverse of dispersion index.

We consider a scaled Poisson approximation to $S$, 
i.e. $S\approx_d \frac{1}{\beta} A_{\beta\mu}$, as introduced in~\eqref{eqn:PoisApprox}.
It can be easily verified that 
the first two moments
of $S$ are matched by that of the approximation. 
We focus on deriving moderate derivations results on the relative error of the tail distributions. 

In order to utilize the Stein-Chen method to compare the two distributions, which only works for integers, 
\textcolor{black}{
we first make the following assumption. We show later our relative error bound also holds without the assumption.}
\begin{asm} \label{asm:1}
Assume that the weights $b_r$'s are positive integers, 
and $\beta$ is rational and can be expressed as $\beta=\frac{n}{m}$ where $n$ and $m$ are relatively prime.
\end{asm}

In some applications such as those in \cite{Sig2012}, Assumption~\ref{asm:1} holds naturally. Under Assumption 1, the scaled Poisson approximation can be rewritten as 
\begin{equation} \label{eqn:nS}
nS \approx_d  \hat{A}_{\lambda}, \quad  \mbox{ where } \hat{A}_{\lambda} := m A_{\lambda},  \mbox{ and }
\lambda = \beta \mu.
\end{equation}
Note that both $nS$ and $\hat{A}_\lambda$ are now integer-valued r.v.s. 

\smallskip
\subsection{Stein's Method for Scaled Poisson Random Variable}
\label{sec:modifiedStein}
We next present a modified Stein's method that will help establish moderate deviations results on the tail distributions under approximation~\eqref{eqn:nS}.
%
%
Note that
the conventional Stein's operator $\mathcal{A}$ given by \eqref{eqn:SteinOp1} 
no longer works
for $\hat{A}_\g$,
since it is a scaled Poisson random variable that only takes values on grid $m \Z$, where $\Z$ denotes the set of nonnegative integers.
The next Lemma establishes the Stein's operator for $\hat{A}_\lambda$. Its proof can be found in Appendix \ref{sec:SteinEqnNew}.
\begin{lemma}\label{lem:SteinEqnNew}
The Stein's operator for $\hat{A}_{\lambda}$ is 
\begin{equation}\label{eqn:stein0}
\mathcal{A}f(w) := \lambda m f(w+m)-wf(w). 
\end{equation}
That is, $E(\mathcal{A}f)(\hat{A}_{\lambda})=0$ for 
all real-valued functions $f$ defined on $Z_{\ge 0}$.
\end{lemma}
%

For a given metric function $h(\cdot)$, the Stein's equation for $\hat{A}_{\lambda}$ bears the following form,
\begin{equation}\label{eqn:steinEqM}
\mathcal{A}f(w) = \lambda m f(w+m)-wf(w)=h(w)-Eh(\hat{A}_{\lambda}).
\end{equation}

The following Lemma provides a representation of $f_h$, the solution to Stein's equation \eqref{eqn:steinEqM},
which will be useful in establishing our moderate deviations results.
The proof can be found in Appendix~\ref{sec:Appendix_A0}.
\begin{lemma} \label{lem:0}
The following form 	
\begin{equation}\label{eqn:stein_soln}
	f_h(w) = -\sum_{j=0}^\infty \frac{(\g m)^j}{\prod_{\ell=0}^j (w+m \ell) } [h(w+mj) -Eh ({\hat A}_\g)]
\end{equation}
is a solution to the Stein's equation~\eqref{eqn:steinEqM}.
\end{lemma}
Based on \eqref{eqn:steinEqM},  the difference $E[h(nS)]-E[h( {\hat A}_\g)]$ can then be estimated, by studying $E\mathcal{A}f_h(nS)= E[\lambda m f(nS+m)-nSf(nS)]$. 
The following lemmas provide keys to such estimations.  
\begin{lemma}
\label{lem:Poisson_sum_decom}
For any bounded function $f(\cdot)$ defined on $\Z$,
\begin{align}
\label{eqn:Poisson_sum_decom}
E[nSf(nS)]=& \sum_{r=1}^R nb_r \nu_rE\left[f\left(nS+nb_r\right)\right]
=m\g \sum_{r=1}^R \delta_r E\left[f\left(nS+nb_r\right)\right],
\end{align}
where $\delta_r = \frac{b_r\nu_r}{\mu}$, and $\sum_{r=1}^R \delta_r = 1$.
\end{lemma}
Proof can be found in Appendix~\ref{sec:appendix_dec}.
The following lemma is then immediate.
\begin{lemma} \label{lem:decomp}
\[E(\mathcal{A}f_h)(nS) =
\g m E[f(nS+ m)] -E[nSf(nS)] 
=H_1+\cdots+H_R,\]
where
\begin{align} 
H_r &= m\g \delta_r E[ (f_h(nS+m ) -f_h(nS+nb_r)], \quad r=1, \dots, R.  \label{eqn:Hr}
\end{align}
\end{lemma}

Based on \eqref{eqn:parameters}, we can write 
\begin{equation} \label{eqn:1beta}
\frac{1}{\beta}=\frac{m}{n}=\frac{\sigma^2}{\mu}=
\sum_{r=1}^R b_r\frac{b_r\nu_r}{\mu}=
\sum_{r=1}^R b_r \delta_r.   
\end{equation}
That is, $m= \sum_{r=1}^R \delta_r\cdot nb_r$, is a weighted sum of $nb_r$'s, with total weights $\delta_r$'s summing to $1$.

%






\smallskip
\subsection{Main Results on Relative Errors for Scaled Poisson Approximation}
\label{sec:moderateResults}
Consider the following notations adapted from 
\cite{liu_xia_2020}. Let
$\pi_k= P(A_{\lambda}=k)= e^{-\lambda} \cdot \frac{\lambda^k}{k!}$, $k=0, 1, \ldots$, and 
\begin{align} \label{eqn:Fn}
F(n) = \sum_{k=0}^n \pi_k \quad \mbox{  and  } \quad {\overline{F}}(n) = \sum_{k=n}^{\infty}\pi_k.    
\end{align}

Using the modified Stein's operator \eqref{eqn:stein0}, we establish the moderate deviations in the scaled Poisson approximation. Our main results are stated in the following theorems.

\begin{theorem}\label{thm:main} 
Consider a weighted Poisson sum 
$S=\sum_{r=1}^R b_rA_{\nu_r}$ with integer weights $0< b_1 < b_2 < ...< b_R$, mean $\mu$ and variance $\sigma^2$.  
Let $\beta=\frac{\mu}{\sigma^2}$ and  $\lambda = \beta \mu$, and 
$\delta_r=\frac{b_r \nu_r}{\mu}$, $r=1, ..., R$.
Assume $\g \ge 1$ and $\beta$ satisfies Assumption~\ref{asm:1}.
Then, 
for all integers $y > \g+1$,  
\begin{align} \label{eqn:main_new}
\Bigg|\frac{P(\beta S\ge y)}{P(A_\g \ge y )}-1\Bigg|
	 \le  (1 + \alpha^*) \Cplus 
 + (1+\alpha^*)\Cminus  +  \alpha^*
\end{align}
where
\begin{align} \label{eqn:defn_C+}
\Cplus &:= \frac{\g+1}{ y\pi_{y-1}} \\
\Cminus &:= 
\frac{2\g}{\min\{ \pi_0, \lfloor \g \rfloor \pi_{y-1}\}}.
	\label{eqn:defn_C-}
\end{align}
and 
\begin{align} \label{eqn:alpha*}
\alpha^* = \sum_{r=r^*+1}^R (\lceil \beta b_r\rceil-2)\delta_r, 
\quad \mbox{with }
r^* := \max\{r: \beta b_{r}\le 1,  r=1, ..., R\}.
\end{align}
\end{theorem}

The proof of Theorem~\ref{thm:main} will be postponed to Sec.~\ref{sec:proofs}. Next, 
we show that the above relative bound also holds when $\beta$ is irrational.

\begin{cor}\label{thm:main_cor} 
Consider a weighted Poisson sum 
$S=\sum_{r=1}^R b_rA_{\nu_r}$ with integer weights $0< b_1 < b_2 < ...< b_R$, mean $\mu$ and variance $\sigma^2$.  
Let $\beta=\frac{\mu}{\sigma^2}$ and  $\lambda = \beta \mu$, and 
$\delta_r=\frac{b_r \nu_r}{\mu}$, $r=1, ..., R$.
Assume $\g \ge 1$ and $\beta \in \R$. Then, \eqref{eqn:main_new} holds for all integers $y > \g+1$.
%
\end{cor}
\begin{proof}
If $\beta$ is rational, then there is nothing to prove. If $\beta$ is irrational, then there is a sequence of rational numbers $\beta_n, n=1,2,\ldots, \infty$, satisfying $\beta_n< \beta$ and $\beta_n\rightarrow \beta$, as $n\rightarrow\infty$. For each $\beta_n$, Theorem \ref{thm:main} implies that, 
\begin{align*}
\Bigg|\frac{P(\beta_n S\ge y)}{P(A_{\g_n} \ge y )}-1\Bigg|
	 \le  (1 + \alpha_n^*) \Cplus 
 + (1+\alpha_n^*)\Cminus  +  \alpha_n^*
\end{align*}
with $\g_n=\beta_n\mu$, and,
\begin{align*} 
\alpha_n^* = \sum_{r=r_n^*+1}^R (\lceil \beta_n b_r\rceil-2)\delta_r, 
\quad \mbox{with }
r_n^* := \max\{r: \beta_n b_{r}\le 1,  r=1, ..., R\}.
\end{align*}
The assumption $\beta_n < \beta$ implies that $r_n^*\ge  r^*$, which leads to $\al^*_n\le \al^*$. Thus,
\begin{align*}
\Bigg|\frac{P(\beta_n S\ge y)}{P(A_{\g_n} \ge y )}-1\Bigg|
	 \le  (1 + \alpha^*) \Cplus 
 + (1+\alpha^*)\Cminus  +  \alpha^*
\end{align*}
It is also easy to see that, as a consequence of the bounded convergence theorem,
\begin{align*}
\lim_{n\rightarrow \infty} P(\beta_n S\ge y)=P(\beta S\ge y), \quad \lim_{n\rightarrow \infty}
P(A_{\g_n} \ge y )=P(A_{\g} \ge y ).
\end{align*}
Inequality~\eqref{eqn:main_new} thus follows.
\end{proof}


The next theorem presents our main result in its final form, where all requirements in Assumption~\ref{asm:1} are removed. That is, our relative error bound also holds for non-integer or irrational $b_r$'s.

\begin{theorem}[{\bf Main Result}]\label{thm:main_cor_II} 
Consider a weighted Poisson sum 
$S=\sum_{r=1}^R b_rA_{\nu_r}$ with weights $0< b_1 < b_2 < ...< b_R$, mean $\mu$ and variance $\sigma^2$.  
Let $\beta=\frac{\mu}{\sigma^2}$ and  $\lambda = \beta \mu$, and 
$\delta_r=\frac{b_r \nu_r}{\mu}$, $r=1, ..., R$.
Assume $\g \ge 1$ and $\beta\in \R$. Then, 
\eqref{eqn:main_new} holds 
for all integers $y > \g+1$.  
\end{theorem}

\begin{proof}
First, consider the case that $b_1, \ldots b_R$ are all  rationals. 
We can find a large enough integer $B$ such that $b_r = \frac{\hat{b}_r}{B}$, and $\hat{b}_r$'s are all integers,  for $r=1, ..., R$.  Let 
\[\hat{\mu} = \sum_{r=1}^R\hat{b}_r\nu_r, \quad \quad \hat{\sigma}^2=\sum_{r=1}^R \hat{b}_r^2\nu_r.\]
We then have,
$\mu =\frac{ \hat{\mu} }{B}, \sigma^2 = \frac{ \hat{\sigma}^2 }{B^2}$, and ${\beta} = B \hat{\beta}$, $ \lambda  = \beta\mu = \hat{\beta}\hat{\mu} =\hat{\lambda}$. 
Therefore, showing  $S \approx_d  \frac{1}{\beta} A_{\beta\mu}$ is equivalent to showing that $ \hat{S} \approx_d  \frac{1}{\hat{\beta}} A_{{\hat{\beta}}{\hat{\mu}}}$, where $\hat{S} = \sum_{r=1}^R \hat{b}_rA_{\nu_r}$ now satisfies that conditions of Corollary 1, 
with all the weights being positive integers.

If some of $b_r$'s are irrational, then there is a sequence of rational $R$-tuples $(b_r^{(n)})_{r=1,\ldots, R}$ such that $b_r^{(n)} \le b_r$, for all $n=1, 2, \ldots, \infty$ and $b_r^{(n)}\rightarrow b_r$, as $n\rightarrow\infty$. For each $n$, Corollary  \ref{thm:main_cor} implies that, 
\begin{align*}
\Bigg|\frac{P(S\ge y)}{P(A_{\g} \ge y )}-1\Bigg|
	 \le  (1 + \alpha_{(n)}^*) \Cplus 
 + (1+\alpha_{(n)}^*)\Cminus  +  \alpha_{(n)}^*
\end{align*}
with
\begin{align*} 
\alpha_{(n)}^* = \sum_{r=r_{(n)}^*+1}^R (\lceil \beta b^{(n)}_r\rceil-2)\delta_r, 
\quad \mbox{with }
r_{(n)}^* := \max\{r: \beta_n b^{(n)}_{r}\le 1,  r=1, ..., R\}.
\end{align*}
The assumption $b_r^{(n)} \le b_r$ implies that $r_{(n)}^*\ge  r^*$, which leads to $\al^*_{(n)}\le \al^*$. 
Inequality~\eqref{eqn:main_new} thus follows.
\end{proof}

\smallskip
\begin{remark}
Theorem~\ref{thm:main_cor_II} indicates that the relative error of the scaled Poisson approximation depends on two main quantities of $\Cplus$ and $\Cminus$. 
They are similar to the quantities ${\mathfrak C}_1(\g,k)$ and ${\mathfrak C}_2(\g,k)$ in \cite{liu_xia_2020} which play crucial roles in bounding the norm of the differences of the solution to the Stein equation for a Poisson approximation.
The numerical results below demonstrate that $\Cplus$ and $\Cminus$ are comparable to, and in many cases smaller than, ${\mathfrak C}_1(\g,k)$ and ${\mathfrak C}_2(\g,k)$. This confirms the power of Stein method in obtaining moderate deviations results for approximation of Poisson type. 
\end{remark}

\medskip

In~\cite{liu_xia_2020}, ${\mathfrak C}_{1}(\g,k)$ and ${\mathfrak C}_{2} (\g,k)$ are compared to a``naive'' counterpart $(1-e^{-\g})/(\g P(A_\g\ge k))$, which bounds the absolute error of the Poisson approximation by $(1-e^{-\g})/\g$ through total variation analysis that can be found in~\cite{46218670-2a65-3e35-9463-f565a234abb9},~\cite{10.1214/aop/1176991491}, and~\cite{BHJBook1992}.
The same comparison can be carried out for $\Cplus$ and $\Cminus$. More specifically, consider two groups of ratios:
\begin{equation}
    \text{ratio } + := \Cplus/[(1-e^{-\g})/(\g P(A_\g\ge y))]
\end{equation}
\begin{equation}
    \text{ratio } - := \Cminus/[(1-e^{-\g})/(\g P(A_\g\ge y))]
\end{equation}
for $y \in \R_{>0}$; and
\begin{equation}
    \text{ratio } i := {\mathfrak C}_{i}(\g,k)/[(1-e^{-\g})/(\g P(A_\g\ge k))], \quad i=1, 2,
\end{equation}
for $k \in \Z$, where ${\mathfrak C}_{1}(\g,k)$ and ${\mathfrak C}_{2} (\g,k)$ were defined in  \cite{liu_xia_2020}, with
\begin{align}
{\mathfrak C}_{1+}(\g, k)=& \frac{F(k-1)}{k\pi_k}\left(1-\frac{{\overline{F}}(k+1)}{{\overline{F}}(k)}\cdot \frac{k}{\g }\right), \label{eqn:C1+}\\
{\mathfrak C}_{1-}(\g, k)=& \frac{F(k-1)}{k\pi_k}\left(1-\frac{F(k-2)}{F(k-1)}\cdot \frac{\g}{k-1}\right), \label{eqn:C1-}
\end{align}
and whose maximum and summation are ${\mathfrak C}_{1}(\g,k)$ and ${\mathfrak C}_{2}(\g,k)$, respectively. 

\begin{figure}[htbp] 
	\centering
	\includegraphics[scale = 0.8]{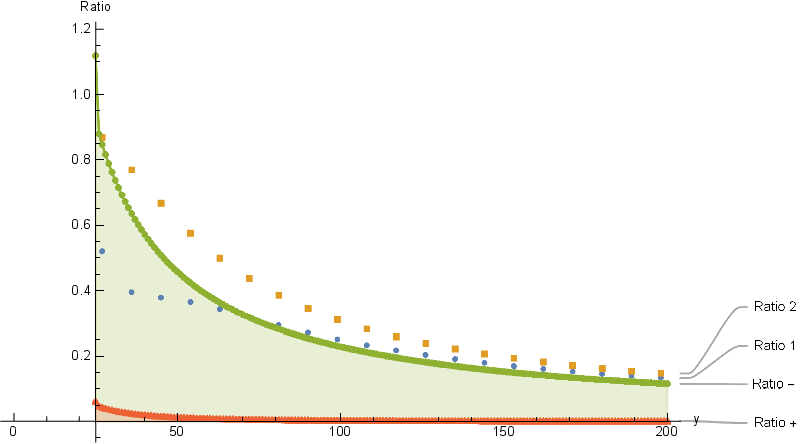}
 	\caption{Comparison of ratio+, ratio-, ratio 1 and ratio 2 for $y \in [25, 200]$ and $k = \lfloor\frac{y}{m}\rfloor$ with $m=9$.}
	\label{fig:ratio}
\end{figure}
\begin{example}
Figure~\ref{fig:ratio} demonstrates the two groups of ratios as a function of the tail parameter~$y$
for the following 2-class example. Consider $S=A_{\nu_1} + 2 A_{\nu_2}$ with $\nu_1 = 1$ and $\nu_2=2$. In this case, $b_1 = 1$, $b_2=2$, and $R=2$. 
Based on \eqref{eqn:parameters}, we have $\mu=5$, $\sigma^2 = 9$, and $\beta=5/9$. Note that
$\beta$ is a rational number with $n=5$ and $m=9$ that are relatively prime. 
Furthermore, we have $\lambda =\beta \mu= 25/9$.

Based on \eqref{eqn:alpha*}, we have $r^*=1$ and $\alpha^*=0$. Thus
the relative error bound \eqref{eqn:main_new} becomes:
\begin{align}
\Bigg|\frac{P(\beta S\ge y)}{P(A_\g \ge y )}-1\Bigg|
	 \le   \Cplus +   \Cminus 
\end{align}

As shown in the figure,  both of our ratios (ratio + and ratio -) perform at similar scale as ratios~1 and~2 in \cite{liu_xia_2020}, relative to the "naive" counterparts.  
Note that ratios $1$ and $2$ are Poisson approximation that defined only on Poisson grid of  $\hat{A}_\lambda(= m A_\lambda)$. Therefore, in the figure, we only have ratios defined for $y$'s that are divisible by $m$, i.e., 
$k = \lfloor\frac{y}{m}\rfloor$ in  \eqref{eqn:C1+} and \eqref{eqn:C1-}.
\end{example}

The following corollary is immediate. It is easily checked from
\eqref{eqn:defn_C+} and \eqref{eqn:defn_C-}. 
\begin{cor}\label{thm:mainII_cor} We have
\begin{equation} \label{eqn:C-Dom}
\Cminus > \Cplus,  \quad 
\mbox{when $y$ is large.}
\end{equation}
That is, the relative error bound \eqref{eqn:main_new} will eventually be dominated by $\Cminus$ when $y$ is large enough. 
\end{cor}

\smallskip
\subsection{On a conjecture in ~\cite{liu_xia_2020}}
\label{sec:conjecture}

In~\cite{liu_xia_2020}, in an effort to establish moderate deviations for (non-scaled) Poisson approximations,  the authors show that the Stein method is reduced to the evaluation of the two quantities, 
${\mathfrak C}_{1-}(\g,k)$ and ${\mathfrak C}_{1+}(\g,k)$, which are defined in \eqref{eqn:C1-} and \eqref{eqn:C1+}, respectively.

In~\cite{liu_xia_2020}, it was conjectured that (Conjecture 4.2) 
\[{\mathfrak C}_{1-}(\g,k)>{\mathfrak C}_{1+}(\g,k), 
\] 
thus ${\mathfrak C}_{1-}(\g,k)$ plays a dominant role in their error bounds for the non-scaled Poisson approximation. 

The next proposition shows that the above conjecture is indeed true when $k$ is large enough and further provides the specific threshold on $k$. The detailed proof can be found in Appendix~\ref{sec:conj}.
 
\begin{proposition}
\label{pro:LX_conjecture}
When $k\ge (\g+1)+\sqrt{(\g+1)(5\g+1)}$, we have,
\[{\mathfrak C}_{1-}(\g,k)>{\mathfrak C}_{1+}(\g,k),
\] 
where ${\mathfrak C}_{1+}(\g,k)$ and ${\mathfrak C}_{1-}(\g,k)$ are defined in  \eqref{eqn:C1+} and 
\eqref{eqn:C1-}, respectively. That is, the moderate deviation results of ~\cite{liu_xia_2020} can be rewritten in terms of   ${\mathfrak C}_{1-}(\g,k)$ and ${\mathfrak C}_{1+}(\g,k)$, and the relative error bound will be eventually dominated by ${\mathfrak C}_{1-}(\g,k)$ for large $k$. 
\end{proposition}

\section{Proof of Theorem~\ref{thm:main}}
\label{sec:proofs}
%
Let 
\begin{equation} \label{eqn:hmy}
 h(w)=\mathbb{I}\{w: w\ge my\}, \quad \mbox{ for } w \in \Z.
\end{equation}
Let $f_h$ be the solution to Stein's equation \eqref{eqn:steinEqM} with $h$ given by \eqref{eqn:hmy}. Since $E[h(nS)]= P(nS \ge my)$, we have
\begin{equation} \label{eqn:Efh}
P(nS \ge my) - P(\hat{A}_\lambda \ge my) 
= E h(nS) -E h({\hat A}_\g) = E\mathcal{A} f_h(nS). 
\end{equation}
Apply Lemma \ref{lem:decomp} to equation \eqref{eqn:Efh}, we have:
\begin{equation} \label{eqn:decompR}
P(nS \ge my) - P(\hat{A}_\lambda \ge my)  = H_1 + \cdots + H_R,
\end{equation}
where $H_r$ is defined in \eqref{eqn:Hr}.
Therefore, to bound the deviation between the tail distributions of $n S$ and $\hat{A}_{\lambda}$, it suffices to bound $|H_1|$, $\cdots$, and $|H_R|$.

The next lemma provides important bounds on differences of $f_h$  which will be useful in bounding $|H_r|$'s.
It depends on a sequence of lemmas 
(Lemmas~\ref{lem:difference}-\ref{lem:Bdminus}) that reveal the key properties of function $f_h$ which are presented in section~\ref{sec:properties}.
The proof of Lemma~\ref{lem:md_all} is provided at the end of section~\ref{sec:properties}.
 
\begin{lemma}\label{lem:md_all} 
For all integers $y \ge \g +1$, we have:
\begin{align}
E[|f_h(nS+m ) -f_h(nS) |]
\le  P(A_\g\ge y)\left[\Bplus +\BMminus \right],   \label{eqn:mgap}
\end{align}
and for $0<\ell<m$,
\begin{align}
& E[|f_h(nS+m)-f_h(nS+\ell)|] 
\le  P(A_\g\ge y) \left[ \Bplus +\Bdminus\right], 
  \label{eqn:dgap}
\end{align}
where 
\begin{align}   
\Bplus & 
=  \frac{\g+1}{\g m y\pi_{y-1}}
\label{eqn:Bplus} \\
\BMminus &= \frac{2}{m \cdot \min\{\pi_0, 
\lfloor \g \rfloor \pi_{y-1}}\}+ \frac{1}{m\g}, \\
\Bdminus &= \frac{2}{m \cdot \min\{ \pi_0, \lfloor \g \rfloor \pi_{y-1}\}}.
\end{align}
\end{lemma}
  
The next lemma presents bounds on the $|H_r|$'s.
\begin{lemma}\label{lem:H_r}
Let 
$r^* = \max\{r: \beta b_{r}\le 1,  r=1, ..., R\}$, and
$K_r = \lceil \beta b_r\rceil$.
\begin{itemize}
\item For $1\le r \le r^*$, we have
\begin{align} \label{eqn:Hr_low}
		|H_r|\le   P(A_\g\ge y) \cdot \delta_r \g m\left[ \Bplus+ \Bdminus
  \right]; 
		\end{align}
\item For $r>r^*$, we have
\begin{align} \label{eqn:Hr_high}
|H_r| \le 
P(A_\g \ge y) \cdot 
\delta_r \g m 
\left[ 
\Big(\Bplus +\Bdminus \Big) + 
(K_r-2)\Big(\Bplus +\BMminus \Big) 
\right].
\end{align}
\end{itemize}
\end{lemma}

\proof  Since $\beta=n/m$, 
 $\{\beta b_{r}\le 1\} \Longleftrightarrow 
 \{nb_r\le m\}$. 
 
 Thus, $r^*= \max\{r: \beta b_{r}\le 1,  r=1, ..., R\}$
is simply the largest $r$ that satisfies $n b_r \le m$; and  for all $r =r^*+1, ..., R$, we have $n b_r > m$ and 
 $K_r=\lceil \beta b_r\rceil \ge 2$.


 For $1 \le r \le r^*$, since $nb_r\le m$, we have $0\le\Delta=m-nb_r<m$. 
 Thus 
 \begin{align*}
|H_r| =  & \g m\delta_rE[|f_h(nS+m)-f_h(nS+nb_r)| \\
  \le  & \g m \delta_r P(A_\g\ge y)  \left[ \Bplus +\Bdminus
  \right],
\end{align*}
where the inequality follows from Lemma~\ref{lem:md_all}. 

 For $r =r^*+1, ..., R$,  since $nb_r>m$, 
 the gap between $nS+m$ and $nS+nb_r$ may contain multiple $m$-steps.  
Since 
$K_r-1 < \beta b_r \le K_r$, and $\beta=\frac{n}{m}$,
 we have
$(K_r-1)m < nb_r \le K_r m$. 
Therefore,
\begin{align*}
|H_r| = 
 & \g m \delta_r E[|f_h(nS+nb_r)-f_h(nS+m)|] \nonumber \\
\le & \g m \delta_r E[|f_h(nS+nb_r)-f_h(nS+(K_r-1)m)|+ \ldots \nonumber \\
& \hspace*{2.75in} + |f_h(nS+2m)-f_h(nS+m)|] \nonumber \\
\le & P(A_\g \ge y) \cdot 
\g m \delta_r 
\left[ 
\Big(\Bplus +\Bdminus \Big) + 
(K_r-2)\Big(\Bplus +\BMminus \Big) 
\right] 
\end{align*}
where the last inequality is due to Lemma~\ref{lem:md_all}.

\endproof


We are now ready to prove Theorem~\ref{thm:main}.

\begin{proof}[\bf  Proof of Theorem~\ref{thm:main}]
Based on \eqref{eqn:decompR} and Lemma~\ref{lem:H_r}, we have
\begin{align*}
&\Big|\frac{P[nS\ge my]}{P[{A}_\g \ge y ]}-1\Big| \\
\le  & \frac{\sum_{r=1}^{r^*} |H_r|+\sum_{r=r^*+1}^R |H_r|}{P(A_\g \ge y)}
\\
\le & 
\sum_{r=1}^{r^*} \delta_r \g m \Big( \Bplus+ \Bdminus\Big)\nonumber\\ & + \sum_{r=r^*+1}^R \delta_r \g m 
\left[ 
\Big(\Bplus +\Bdminus \Big) + 
(K_r-2)\Big(\Bplus +\BMminus \Big) 
\right]  \nonumber \\
= & \Big(1 + \sum_{r^*+1}^R (K_r-2)\delta_r\Big) \g m\Bplus) +\g m \Bdminus
+ \Big(\sum_{r^*+1}^R (K_r-2)\delta_r\Big) \g m \BMminus \\
= & (1 + \alpha^*) \g m  \frac{\g+1}{\g m y\pi_{y-1}}
+ \frac{2\g m }{m \cdot \min\{ \pi_0,  \lfloor \g \rfloor \pi_{y-1}\}}
  +  \alpha^*\g m \Big[\frac{2}{m \cdot \min\{\pi_0, 
\lfloor \g \rfloor \pi_{y-1}}\}+ \frac{1}{m\g}
  \Big] \\
  = & (1 + \alpha^*)  \frac{\g+1}{ y\pi_{y-1}}
 + (1+ \alpha^*)\frac{2\g}{ 
 \min\{ \pi_0,  \lfloor \g \rfloor \pi_{y-1}\}}
 +  \alpha^* \\
  = & (1 + \alpha^*) \Cplus) 
 + (1+ \alpha^*)\Cminus  +  \alpha^*
    \end{align*}
where $\alpha^*$ is defined by \eqref{eqn:alpha*}, and
$\Cplus$ and $\Cminus$ by \eqref{eqn:defn_C+} and \eqref{eqn:defn_C-} respectively.
Theorem~\ref{thm:main} therefore follows.
\end{proof}


\section{Properties of the Solution to the Stein Equation}
\label{sec:properties}


In this section, we present a sequence of lemmas that provide important properties of function $f_h$,  the solution to the Stein's equation, which are essential to the proof of Lemma~\ref{lem:md_all}.

\subsection{Some Basic Properties of $f_h$} \label{subsec:basic} 
Recall that for $h(w)=\mathbb{I}\{w: w\ge my\}$, where $w\in\Z$, Lemma~\ref{lem:0} states that $f_h$, the solution to the Stein's equation, takes the following form: 
\begin{align} \label{eqn:fhbase}
f_h(w) = -\sum_{j=0}^\infty \frac{(\g m)^j}{\prod_{\ell=0}^j (w+m \ell) } [h(w+mj)-P({\hat A}_\g \ge my)]. 
\end{align}

Define $k_w:=\lceil \frac{(my-w)^+}{m}\rceil$. Equivalently, $k_w=\inf\{k\in\Z, w+km \ge my\}$, i.e. $k_w$ is the minimum number of $m$-jumps, starting from $w$, to reach or go above $my$ (i.e. $\ge my$). When $w \ge my$, $k_w=0$ naturally. 
Therefore, $h(w+mj)=\mathbb{I}\{j\ge k_w\}$.
Moreover, $P({\hat A}_\g \ge my)=P({A}_\g \ge y)$. Thus \eqref{eqn:fhbase} can be rewritten as, 
\begin{align}  \label{eqn:SteinEq3}
f_h(w) =&P({A}_\g \ge y) \sum_{j=0}^{\infty}\frac{(\g m)^j}{\prod_{\ell=0}^j (w+m \ell) }  -\sum_{j=k_w}^\infty \frac{(\g m)^j}{\prod_{\ell=0}^j (w+m \ell) }.
\end{align}
Denote
\begin{align} \label{eqn:Thetas}
\Theta(k, w):=\sum_{j=0}^{k} \frac{(\g m)^j}{\prod_{\ell=0}^j(w+ m\ell) },
\qquad 
\overline{\Theta}(k,w):=\sum_{j=k}^\infty \frac{(\g m)^j}{\prod_{\ell=0}^j(w+ m\ell) }, \qquad \mbox{ for $k\ge 0$}
\end{align}
\eqref{eqn:SteinEq3} can then be rewritten as follows:
\begin{align}  
f_h(w) =&P({A}_\g \ge y) \overline{\Theta}(0,w) -  \overline{\Theta}(k_w, w) 
\label{eqn:SteinEq5} 
\end{align}
The next lemma provides an important relation between $f_h(w+m)$ and $f_h(w)$.
\begin{lemma}
\label{lem:difference} 
For all $w \in \Z$,
\begin{align}
\label{eqn:difference_rel_I}
f_h(w+m)= \frac{w}{\g m} f_h(w) - \frac{P({A}_\g\ge y) -1_{w \ge my}}{\g m} 
\end{align}
or, equivalently,
\begin{align}
\label{eqn:difference_rel_II}
Df_h(w):=f_h(w+m) -f_h(w)  = &f_h(w)  \left(\frac{w}{\g m}-1\right) - \frac{P({A}_\g\ge y)-1_{w \ge my}}{\g m}.
\end{align}

\end{lemma}
\proof 
If $w \ge my$, then $k_w=0$ and $k_{w+m}=0$. Then, \eqref{eqn:SteinEq3} becomes
\begin{align*}
f_h(w+m) &=- P({A}_\g < y)\sum_{j=0}^\infty\frac{(\g m)^j}{\prod_{\ell=0}^j (w+m (\ell+1)) } \\
&=- P({A}_\g < y)\sum_{j=0}^\infty\frac{(\g m)^{j+1}}{\prod_{\ell=0}^{j+1} (w+m \ell) }
\cdot \frac{w}{\g m} \\
&=- P({A}_\g < y)\left[\sum_{j'=0}^\infty\frac{(\g m)^{j'}}{\prod_{\ell=0}^{j'} (w+m \ell) }- \frac{1}{w}\right]\cdot \frac{w}{\g m}\\
&= \frac{w}{\g m} f_h(w) + \frac{P({A}_\g < y)}{\g m}
\end{align*}

If $w < my$, then $k_w=\lceil \frac{my-w}{m}\rceil\ge 1$, and $k_{w+m}=k_w -1$.  Then \eqref{eqn:SteinEq3} becomes
\begin{align*}   
&f_h(w+m)  \nonumber \\ 
= & P({A}_\g\ge y) \sum_{j=0}^{\infty}\frac{(\g m)^j}{\prod_{\ell=0}^j (w+m+m \ell) }- \sum_{j=k_w-1}^\infty\frac{(\g m)^j}{\prod_{\ell=0}^j (w+m+m \ell) } \nonumber \\
=& P({A}_\g\ge y) \sum_{j=0}^{\infty}\frac{(\g m)^{j+1}}{\prod_{\ell=0}^{j+1} (w+m \ell) }\cdot \frac{w}{\g m} - \sum_{j=k_w-1}^\infty\frac{(\g m)^{j+1}}{\prod_{\ell=0}^{j+1} (w+m \ell) } \cdot \frac{w}{\g m} \nonumber \\
=&P({A}_\g\ge y) \left[\sum_{j=0}^{\infty}\frac{(\g m)^{j}}{\prod_{\ell=0}^{j} (w+m \ell) } -\frac{1}{w}\right] \cdot \frac{w}{\g m} - \sum_{j=k_w}^\infty\frac{(\g m)^{j}}{\prod_{\ell=0}^{j} (w+m \ell) } \cdot \frac{w}{\g m} \nonumber \\ 
=&f_h(w)  \left(\frac{w}{\g m}\right) - \frac{P({A}_\g\ge y)}{\g m}.
\end{align*}
In both cases, equation \eqref{eqn:difference_rel_I} is satisfied.\\
\endproof

The following lemma will be useful in deriving bounds. The proof is in Appendix~\ref{sec:Lemma14}. 

\begin{lemma}\label{lem:quantity_two} 
For any integer $k\ge 0$, 
we have:
\begin{align} \label{eqn:barTheta}
\frac{1}{m\g \pi_{\lfloor\frac{w}{m}\rfloor}}
P\left[A_\g \ge k+\left\lfloor\frac{w}{m}\right\rfloor+1\right]
\le { \overline{\Theta}(k, w)}
&\le
\frac{1}{\max(1, \lfloor \frac{w}{m} \rfloor) \cdot 
m \pi_{\lfloor\frac{w}{m}\rfloor}}
P\left[A_\g \ge k+\left\lfloor\frac{w}{m}\right\rfloor\right].
\end{align}
and
\begin{align} \label{eqn:Theta}
\frac{1}{m\g \pi_{\lfloor\frac{w}{m}\rfloor}}
P\left(\left\lfloor\frac{w}{m}\right\rfloor+1
\le A_\g \le k+\left\lfloor\frac{w}{m}\right\rfloor+1\right)
\le 
{{\Theta}(k, w)}
& \le 
\frac{1}{\max(1, \lfloor \frac{w}{m} \rfloor) \cdot 
m\pi_{\lfloor\frac{w}{m}\rfloor}}
P\left(\left\lfloor\frac{w}{m}\right\rfloor\le 
A_\g \le k+\left\lfloor\frac{w}{m}\right\rfloor\right).
\end{align}
\end{lemma}

The following two lemmas provide some useful results on Poisson probabilities. 
\begin{lemma} \label{lem:Pois}
For any Poisson r.v. $A_\g$ with $\g >0$, we have
\begin{align}
 \frac{P(A_\g \ge k+1)}{P(A_\g \textcolor{black}{\ge} k)} 
& \le \frac{\g}{k+1},   \mbox{ for any $k > \g-1$}  \label{eqn:PoisB1}\\
\frac{P(A_\g \ge y+k)}{P(A_\g \ge y)} 
& \le \frac{\g}{y+k},  \mbox{ for any $k\ge 1$, and $y\ge \g$}  \label{eqn:PoisB2}
\end{align}
\end{lemma}

\proof
It has been observed in \cite{ChenFangShao2013} that 
\begin{align*}
P(A_\g \ge k) =& P(A_\g = k) ( 1+ \frac{\g}{k+1} + \frac{\g^2}{(k+1)(k+2)} + ...) \\
\le& P(A_\g = k) ( 1+ \frac{\g}{k+1} + \frac{\g^2}{(k+1)^2} + ...) \\
= &  P(A_\g = k) \frac{k+1}{k-\g + 1}   \quad \mbox{ for all $k > \g -1$.}
\end{align*}

Thus, \eqref {eqn:PoisB1} follows noting that
\begin{align*}
& \frac{P(A_\g \ge k+1)}{P(A_\g \ge k)} 
= 1-\frac{P(A_\g = k)}{P(A_\g \ge k)} 
\le 1 - \frac{k-\g +1}{k+ 1}
= \frac{\g}{k+1}.
\end{align*}

To show \eqref{eqn:PoisB2},  we apply \eqref{eqn:PoisB1} to $y+k-1$, for any $k \ge 1$ and $y\ge \g$, 
\begin{align*}
P(A_\g \ge y+k)\le  \frac{\g}{y+k} P(A_\g \ge y+k-1) 
\le \frac{\g}{y+k} P(A_\g \ge y) 
\end{align*}
\endproof

It is well known that the mode of Poisson distribution is $\lfloor \g \rfloor$ for non-integer $\g$. When $\g$ is a positive integer, the modes are $\g$ and $\g -1$. Furthermore, the Poisson probabilities are increasing before the mean value $\g$ and decreasing after $\g$.
\begin{lemma} \label{lem:Pois2}
For any Poisson r.v. $A_\g$ with $\g >0$ and $\pi_k = P(A_\g = k)$, we have
\begin{align}
 \frac{\pi_k }{\pi_{k-1}}& =\frac{\g}{k}, \quad \quad  \mbox{for any $k\ge 1$; }  \label{eqn:Pik},  \\
 \pi_{\lfloor \frac{w}{m} \rfloor-1} & \ge
\frac{\lfloor \g \rfloor}{\g}\pi_{y-2}, 
\quad \mbox{for $ m\g \le w<my$, where $y$ is an integer and $y>\g + 1$.}
\label{eqn:Piy_2}
\end{align}
\end{lemma}
\proof 
\eqref{eqn:Pik} can be easily verified using the definition of $\pi_k$.

If $\g <1$, \eqref{eqn:Piy_2} is immediate since $\lfloor \g \rfloor=0$.

Now suppose $\g \ge 1$. Since $m\g \le w<my$, we have  
$\lfloor \g \rfloor \le \lfloor \frac{w}{m} \rfloor  \le  y-1$ and
$\lfloor \g \rfloor -1 \le \lfloor \frac{w}{m} \rfloor -1 \le  y-2$.

If $\lfloor \frac{w}{m} \rfloor -1 \ge \lfloor \g \rfloor$, then 
\eqref{eqn:Piy_2} is immediate since $\lfloor \g \rfloor$ is the mode and the Poisson probabilities are decreasing after the mean value $\g$.

If $\lfloor \frac{w}{m} \rfloor -1 = \lfloor \g \rfloor-1$, then
\[
\pi_{\lfloor \frac{w}{m} \rfloor-1} 
= \pi_{\lfloor \g \rfloor-1} 
= \pi_{\lfloor \g \rfloor} \frac{\lfloor \g \rfloor}{\g}
\ge
\frac{\lfloor \g \rfloor}{\g}\pi_{y-2}, 
\]
where the second equality is due to \eqref{eqn:Pik}, and the last inequality is because $y-2\ge \lfloor \g \rfloor$ since 
$y-1>\g$ and $y$ is an integer.

\endproof

\medskip
We next analyze the behavior of $f_h(w)$ for the cases of $w \ge my$ and $w < my$ separately.
\subsection{The Case of $w \ge my$}

For $w\ge my$, observe that $h(w+mj)=1$ for all integer $j \ge 0$, 
and $P({\hat A}_\g < my)= P({A}_\g < y)$,
thus, 
\begin{align}
\label{eqn:solution_for_large_w}f_h(w) = -P({A}_\g < y) \sum_{j=0}^\infty \frac{(\g m)^j}{\prod_{\ell=0}^j (w+m \ell) },  \mbox{ for } w \ge my.
\end{align}
The following lemma is immediate. 
\begin{lemma}\label{lem:mono}
For $w\ge my$, $f_h(w) < 0$ and $f_h(w)$ is monotonically increasing and concave 
in $w$.
\end{lemma}
\proof 
It is easy to verify that each summand in \eqref{eqn:solution_for_large_w}
\begin{align*}
\frac{(\g m)^j}{\prod_{\ell=0}^j (w+m \ell) }
\end{align*}
is decreasing and convex with respect to $w$. 
\endproof
\begin{remark}
Lemma \ref{lem:mono} tells us that the differences of the Stein function, the key quantities in our study, are non-negative and decreasing in $w$. Hence, to estimate the maximum value they can achieve, only the first of them needs to be evaluated, for $w \ge my$.
\end{remark}

\begin{lemma} \label{lem:C+}
For $w\ge my$, we have, 
\[
0 \le \frac{f_h(w+l)-f_h(w)}{l} \le  P(A_\g \ge y) \cdot 
\frac{\g+1}{\g m^2 y\pi_{y-1}}
\quad
\mbox{ for all $l>0$}.
\]
\end{lemma}
\proof 
Since $f_h(w)$ is increasing and concave for $w \ge my$, we have, for all $l>0$, 
\begin{align} \label{eqn:Lem16Eq1}
\frac{f_h(w+l)-f_h(w)}{l} \le \frac{d}{dw} f_h(w)\Big|_{w=my} = P({A}_\g < y)\sum_{j=0}^\infty \frac{(\g m)^j}{\prod_{\ell=0}^j (my+m \ell) }\sum_{k=0}^j \frac{1}{my+mk}.
\end{align}
where the RHS is because
\[\frac{d}{dw} \frac{(\g m)^j}{\prod_{\ell=0}^j (w+m \ell) }
= - \frac{(\g m)^j}{\prod_{\ell=0}^j (w+m \ell) } \sum_{k=0}^j \frac{1}{w+mk}.
\]
Note that 
\begin{align}
&\sum_{j=0}^\infty\sum_{k=0}^j \frac{(\g m)^j}{\prod_{\ell=0}^j (my+m \ell) } \frac{1}{(my+mk)} \nonumber \\ 
= & \sum_{k=0}^\infty\sum_{j=k}^\infty \frac{(\g m)^j}{\prod_{\ell=0}^j (my+m \ell) } \frac{1}{(my+mk)}\nonumber \\
= & \sum_{k=0}^\infty\frac{1}{my+mk}\sum_{j=k}^\infty \frac{(\g m)^j}{\prod_{\ell=0}^j (my+m \ell) } \label{eqn:tight1} \\
= & \frac{1}{m^2}\sum_{k=0}^\infty\frac{1}{y+k}\sum_{j=k}^\infty \frac{\g^j}{\prod_{\ell=0}^j (y+\ell) } \nonumber\\
= & \frac{1}{m^2}\sum_{k=0}^\infty\frac{1}{y+k}\sum_{j=k}^\infty \frac{e^{-\g}\g^{j+ y}}{(y+j)!} \cdot \frac{( y-1)!}{e^{-\g}\g^{y}}  \nonumber \\
= & \frac{1}{\g m^2 \pi_{y-1}}\sum_{k=0}^\infty\frac{1}{y+k}P(A_\g \ge  y  +k)
\label{eqn:tight2}
\end{align}
Based on \eqref{eqn:PoisB2}, we have 
\begin{align*}
\sum_{k=0}^\infty\frac{1}{y+k}\frac{P(A_\g \ge y+k)}{P(A_\g \ge y)}
\le \frac{1}{y} + 
\sum_{k=1}^\infty\frac{1}{y+k}\cdot \frac{\g}{y+k} 
\le \frac{1}{y} + \g(\sum_{k=1}^\infty\frac{1}{y+k} - \frac{1}{y+k+1} ) \le \frac{\g +1}{y}.
\end{align*}
Combine with \eqref{eqn:Lem16Eq1}, 
we then have  
\begin{align*}
\frac{f_h(w+l)-f_h(w)}{l \cdot P(A_\g \ge y)} \le 
\frac{P({A}_\g < y)}{\g m^2 \pi_{y-1}} \cdot \frac{\g+1}{y}
\le 
\frac{\g+1}{\g m^2 y\pi_{y-1}}
\end{align*}

\endproof

\begin{remark}
Note that \eqref{eqn:tight1} can be written as:
\begin{align}
& \sum_{k=0}^\infty\frac{1}{my+mk} \overline{\Theta}(k, my) \nonumber \\
\le & \sum_{k=0}^\infty\frac{1}{m(y+k)} \cdot 
\frac{1}{y \cdot m \pi_{y}}
P\left[A_\g \ge y+ k\right]  \quad \mbox{(based on \eqref{eqn:barTheta} with $w=my$ and $\lfloor \frac{w}{m}\rfloor =y$)} \nonumber \\
= & \sum_{k=0}^\infty\frac{1}{m(y+k)} \cdot 
\frac{1}{m \g \pi_{y-1}}
P\left[A_\g \ge y+ k\right]  \quad \mbox{(since $y\pi_y = \g \pi_{y-1}$)} \nonumber \\
= & \frac{1}{\g m^2 \pi_{y-1}}
\sum_{k=0}^\infty\frac{1}{y+k} \cdot 
P\left[A_\g \ge y+ k\right]  \label{eqn:tight3}
\end{align}
Note that \eqref{eqn:tight3} is identical as \eqref{eqn:tight2}, which indicates that the upper bound given by $\eqref{eqn:barTheta}$ is tight for $w =my$.
\end{remark}

\subsection{The Case of $w< my$}
When $w< my$, $k_w \ge 1$.  Note that we can rewrite \eqref{eqn:SteinEq3} as:
\begin{align}
f_h(w) 
&= P({A}_\g\ge y) \sum_{j=0}^{k_w-1}\frac{(\g m)^j}{\prod_{\ell=0}^j (w+m \ell) }- P({A}_\g < y)\sum_{j=k_w}^\infty\frac{(\g m)^j}{\prod_{\ell=0}^j (w+m \ell) }, 
\quad \forall k_w \ge 1.\label{eqn:split}
\end{align}
Or equivalently, 
\begin{align}
f_h(w) =&P({A}_\g \ge y) \Theta(k_w-1, w) - P({A}_\g < y) \overline{\Theta}(k_w, w), \label{eqn:SteinEq4}
\end{align}
\begin{lemma}
\label{lem:decreasing}
For any $m \lambda \le w < my$, $Df_h(w)=f_h(w+m) -f_h(w)<0$. That is, $f_h(w)$ is decreasing on the $m$-grid.
\end{lemma}
\proof 
If $k_w =1$, we have $w+m \ge my$ and $w<my$. 
Based on Lemma~\ref{lem:difference}, from \eqref{eqn:difference_rel_I},
we have $f_h(w) = \frac{\g m}{w} f_h(w+m)+ \frac{P({A}_\g\ge y)}{w}$.
Thus,
\begin{align*}
f_h(w+m) - f_h(w) =(1-\frac{\g m}{w})f(w+m) - \frac{P({A}_\g \ge y)}{w} < 0,
\end{align*}
where the inequality is because $f(w+m) < 0$ based on Lemma~\ref{lem:mono}.

Now consider the case $k_w\ge 2$. Based on \eqref{eqn:split} and $k_{w+m}=k_w -1$, it is easy to verify that, 
\begin{align*}
f_h(w+m)- f_h(w) =& P({A}_\g \ge y) \sum_{j=0}^{k_w-2}\left[\frac{(\g m)^j}{\prod_{\ell=0}^j (w+m+m \ell)}-\frac{(\g m)^j}{\prod_{\ell=0}^j (w+m \ell)}\right]
\\& - P({A}_\g \ge y) \frac{(\g m)^{k_w-1}}{\prod_{\ell=0}^{k_w-1} (w+m \ell) }
\\ & -(\frac{w}{\g m}-1) \sum_{j=k_w}^\infty \frac{(\g m)^j}{\prod_{\ell=0}^j (w+m \ell) } [1-P({A}_\g \ge y)]  <0,
\end{align*}
where the first two terms are obviously negative, and the last term is negative 
since $m\g \le w < my$.
\endproof

\begin{lemma} 
\label{lem:grid_cocavity}
For $ m\g \le w<my$, 
%
for any $j=1, ..., k_w-1$, where $k_w=\lceil \frac{(my-w)}{m}\rceil$, we have
\begin{align} \label{eqn:2nd_order_diff_k}
D_h(w+jm)-D_h(w+(j-1)m) \le 
P({A}_\g \ge y) \cdot T_m(y),
\end{align}
where 
$T_m(y)=\frac{1}{ m \g
{\lfloor \g \rfloor}\pi_{y-2}}$.
\end{lemma}

\proof 


Note that $k_w=\lceil \frac{(my-w)}{m}\rceil$ 
is the minimum number of $m$-jumps required to reach or go above $my$ starting from $w$.
Hence, $k_w \ge 1$, $w+(k_w-1)m<my$, and $w+k_w m \ge my$. 
For any $j=1,..., k_w-1$, we have $w \le w+(j-1)m < w + jm <my$.

Apply \eqref{eqn:difference_rel_II} of Lemma~\ref{lem:difference} to both $w+(j-1)m$ and $w+jm$, we have, 
\begin{align*}
& Df_h(w+jm)]-Df_h(w+(j-1)m) \nonumber\\
= &f_h(w+jm)  \left(\frac{w+jm}{\g m}-1\right) -f_h(w+(j-1)m)  \left(\frac{w+(j-1)m}{\g m}-1\right) \nonumber \\
=&[f_h(w+jm) -f_h(w+(j-1)m)]\left(\frac{w+jm}{\g m}-1\right) + \frac{f_h(w+(j-1)m)}{\g} \nonumber \\
\le & \frac{f_h(w)}{\g} \quad \quad 
\mbox{(since $f_h(w+jm)$ is decreasing on the m-grid)} \nonumber \\
\le & \frac{P({A}_\g \ge y)\cdot \overline{\Theta}(0, w)}{\g} 
\quad \mbox{(based on \eqref{eqn:SteinEq5})}\nonumber\\  
\le & \frac{P({A}_\g \ge y)}{\g} \frac{P\left[A_\g \ge \left\lfloor\frac{w}{m}\right\rfloor\right]}{
\lfloor \frac{w}{m} \rfloor 
m \pi_{\lfloor\frac{w}{m}\rfloor}} 
\quad \mbox{(based on \eqref{eqn:barTheta}) and 
$1 \le \lfloor \g \rfloor \le \lfloor \frac{w}{m} \rfloor$)}\nonumber\\ 
= & P({A}_\g \ge y) \frac{1}{ m \g^2 \pi_{\lfloor \frac{w}{m} \rfloor-1}}
\quad \mbox{(since $k\pi_k = \g \pi_{k-1}$)} \nonumber\\ 
\le & P({A}_\g \ge y) \frac{1}{ m \g
{\lfloor \g \rfloor}\pi_{y-2}} 
\end{align*}
where the second to last inequality is due to \eqref{eqn:barTheta}, and 
the last inequality is due to \eqref{eqn:Piy_2}. 
\endproof

\medskip
The following three lemmas are important for the proof of the main theorem.
The proofs are available in the Appendix. 
\begin{lemma} \label{lem:smallW}
For $w< \g m$,
\begin{align}
\frac{| f_h(w) -f_h(w+m)| }{P({A}_\g\ge y) }
\le \frac{1}{m \pi_0} + \frac{1}{m\g}.
\label{eqn:w_small}
\end{align}
\end{lemma}
The proof of Lemma \ref{lem:smallW} can be found in Appendix \ref{sec:smallW}.

\begin{lemma}	\label{lem:mediumW}
	For $ m\g \le w < my$,
 \begin{align}
	\label{eqn:g1_bound_new}
 0 \le \frac{f_h(w) -f_h(w+m)}{P({A}_\g\ge y)} \le 
\frac{2}{m\lfloor\g \rfloor\pi_{ y- 1}}  + \frac{1}{my}(\frac{\g+1}{\g}).
\end{align}

\end{lemma}
The proof of Lemma \ref{lem:mediumW} can be found in Appendix~\ref{sec:BMminus}.

Combine the above two lemmas we then have the following result.
\begin{lemma} \label{lem:<my}
For $w<my$, 
	\label{lem:3}
	\begin{align} 
	\frac{|f_h(w+m) -f_h(w)|}{P({A}_\g\ge y)} \le  
\frac{2}{m \cdot \min\{\pi_0, 
\lfloor \g \rfloor \pi_{y-1}}\}+ \frac{1}{m\g}
  =:\BMminus
  \label{eqn:C<my}
	\end{align} 
\end{lemma}

\begin{lemma} \label{lem:Bdminus}
For $w<my$, we have, 
	\begin{align} 
	\frac{|f_h(w+d) -f_h(w)|}{P({A}_\g\ge y)} \le  
 \frac{2}{m \cdot \min\{ \pi_0,  \lfloor \g \rfloor \pi_{y-1}\}}
  =:\Bdminus, 
  \quad \mbox{for $d=1,2,\ldots, m-1$.}
  \label{eqn:Bdw}
	\end{align} 
\end{lemma}
The proof of Lemma \ref{lem:Bdminus} can be found in Appendix \ref{sec:Bdminus}.

\smallskip
We are now ready to prove Lemma~\ref{lem:md_all}.
\proof   {\bf{(Proof of Lemma~\ref{lem:md_all})} }
Based on Lemma~\ref{lem:C+}, we know that, for any $l=0, 1, ..., m-1$, 
\begin{align} \label{eqn:B+}
\frac{E[|f_h(nS+m)-f_h(nS+\ell)| {\bf 1} _{nS+\ell \ge my}]}{P(A_\g \ge y)}\le 
(m-l) \cdot 
\frac{\g+1}{\g m^2 y\pi_{y-1}}
\le
\frac{\g+1}{\g m y\pi_{y-1}}
=:\Bplus.
\end{align}
\eqref{eqn:mgap} follows noting that
\begin{align*}
&E[|f_h(nS+m ) -f_h(nS) |
] \\ 
 = & E[|f_h(nS+m ) -f_h(nS) |\cdot {\bf 1} _{nS \ge my}] +  E[|f_h(nS+m ) -f_h(nS) |\cdot {\bf 1} _{nS < my}] \\
\le & P(A_\g \ge y) \cdot { \Bplus }
+ P(A_\g\ge y)\cdot \BMminus.
\end{align*}
where the inequality is due to \eqref{eqn:B+} and Lemma~\ref{lem:<my}.

\eqref{eqn:dgap} follows because
\begin{align*}
 & E[|f_h(nS+m)-f_h(nS+\ell)|
 ]\\
  = & E[|f_h(nS+m)-f_h(nS+\ell)|({\bf 1} _{nS+\ell \ge my} ]+  E[|f_h(nS+m)-f_h(nS+\ell)|{\bf 1} _{nS+\ell < my})]\\
  \le  &  P(A_\g\ge y) \left[\Bplus +
  \Bdminus\right],
\end{align*}
where the inequality is due to \eqref{eqn:B+} and 
Lemma~\ref{lem:Bdminus}.
\endproof

\section{Conclusions}
\label{sec:conclusions}

In this paper, 
we derive bounds on the relative error of a scaled Poisson approximation to the weighted sum of independent Poisson random variables via a modified Stein's method. While our proof techniques require the studies of a new Stein's operator that acts on the grid $mZ_{\ge 0}$ when $\beta$ is a rational number $\beta=\frac{n}{m}$, our final result does not rely on the rational  assumption. 
\textcolor{black}{
We show that the relative error of the approximation is essentially bounded by two  main quantities $\Cplus$ and $\Cminus$. 
}
The ratios of these two quantities relative to the naive bound perform at similar scale as the ratios in \cite{liu_xia_2020} in the non-scaled Poisson approximation. Furthermore, we demonstrate that one of the conjectures raised in~\cite{liu_xia_2020} is valid beyond a threshold that has an explicitly expression.



\bibliographystyle{plain}

%
%
%

\medskip
\begin{appendix}

\medskip
\section{Proof of Lemma \ref{lem:SteinEqnNew}}
\label{sec:SteinEqnNew}
\proof{Proof} From the definition of $\hat{A}_{\lambda}$, we can see easily that,
\begin{align*}
E(\mathcal{A}f)(\hat{A}_{\lambda})&=\lambda m E[f(m A_{\lambda}+m))-mE[A_{\lambda}f(mA_{\lambda})]]\\
&=\lambda m E[f(m(A_{\lambda}+1))]-m\sum_{j=0}^\infty f(mj)jP(A_{\lambda}=j)\\
&=\lambda m E[f(m(A_{\lambda}+1))]-m\sum_{j=1}^\infty f(mj)\lambda P(A_{\lambda}=j-1)\\
&=\lambda m E[f(m(A_{\lambda}+1))]-\lambda m \sum_{j=0}^\infty f(m(j+1))P(A_{\lambda}=j)\\
&=\lambda m E[f(m(A_{\lambda}+1))]-\lambda m E[f(m(A_{\lambda}+1))]\\
&=0,
\end{align*}
where the third equation makes use of the fact of $\g P(A_{\g}=j-1) =jP(A_\g =j), \forall j \ge 1$.
\endproof

\section{Proof of Lemma \ref{lem:0}}
\label{sec:Appendix_A0}
\noindent
From \eqref{eqn:stein_soln}, we have,
\begin{align*}
\g m f_h(w+m) & = -\g m \sum_{j=0}^\infty \frac{(\g m)^j}{\prod_{\ell=0}^j (w+m +m \ell) } [h(w+m+mj) -E h ({\hat A}_\g)] \\ &= -\sum_{j=0}^\infty \frac{(\g m)^{j+1}}{\prod_{\ell=0}^j (w+m (\ell+1)) } [h(w+m(j+1)) -E h ({\hat A}_\g)]\\
 &= -\sum_{j=0}^\infty \frac{(\g m)^{j+1}}{\prod_{\ell=1}^{j+1} (w+m\ell) } [h(w+m(j+1)) -E h ({\hat A}_\g)]\\ 
 &=-w\sum_{j=0}^\infty \frac{(\g m)^{j+1}}{\prod_{\ell=0}^{j+1} (w+m\ell) } [h(w+m(j+1)) -E h ({\hat A}_\g)]\\ 
 &= -w\sum_{j=1}^\infty \frac{(\g m)^{j}}{\prod_{\ell=0}^{j} (w+m\ell) } [h(w+mj) -E h ({\hat A}_\g)] \\ 
 &=wf_h(w) +  [h(w) -E h ({\hat A}_\g)].
\end{align*}
which satisfies the Stein's equation \eqref{eqn:steinEqM}.

\section{Proof of Lemma \ref{lem:Poisson_sum_decom}}
\label{sec:appendix_dec}
Recall that $S=\sum_{r=1}^R b_r A_{\nu_r}$, therefore,
\begin{align*}
E[nSf(nS)]=& \sum_{r=1}^R nb_r E\left[A_{\nu_r}f\left(\sum_{r=1}^R nb_r A_{\nu_r}\right)\right]\\=&\sum_{r=1}^R nb_r E\left[E\left[A_{\nu_r}f\left(\sum_{r=1}^R nb_r A_{\nu_r}\right)\bigg|A_{\hat \nu_r}\right]\right]
\\=&\sum_{r=1}^R nb_r\nu_r E\left[f\left(\sum_{r=1}^R nb_r A_{\nu_r} +nb_r\right)\right]
\\=&\sum_{r=1}^R nb_r \nu_rE\left[f\left(nS+nb_r\right)\right] 
\\=& m\g \sum_{r=1}^R \delta_r E\left[f\left(nS+nb_r\right)\right] 
\end{align*}
with $A_{\hat \nu_r}:= \{A_{\nu_1}, A_{\nu_2}, \ldots, A_{\nu_R}\} -\{A_{\nu_r}\}$, for $r=1, 2, \ldots, R$. The third equality is based on the fact that the Poisson distribution is in the kernel of operator $\mathcal{A}$ in \eqref{eqn:SteinOp1}, thus for any
Poisson random variable $A_v$ with mean~$v$, we have $E[A_v f(A_v)]= v E[f(A_v+1)]$ for any bounded function $f$; therefore, for each $r=1,\dots, R$, we have 
\begin{align*}
E\left[A_{\nu_r}f\left(\sum_{r=1}^R nb_r A_{\nu_r}\right)\bigg|A_{\hat \nu_r}\right]=&E\left[A_{\nu_r}g_r\left( A_{\nu_r}\right)\bigg|A_{\hat \nu_r}\right]
\\=&\nu_r E\left[ g_r(A_{\nu_r}+1) \right]
\\=&\nu_r E\left[f\left(\sum_{r=1}^R nb_r A_{\nu_r}+nb_r\right)\right],
\end{align*}
with $g_r(w)=f(\sum_{r'=1,r'\neq r}^R nb_{r'} A_{\nu_{r'}}+nb_rw)$.
The last equality is because $m\g = n \mu$, and $\delta_r =\frac{b_r\nu_r}{\mu}$.

\smallskip
\section{Proof of Lemma~\ref{lem:quantity_two} } \label{sec:Lemma14}
\proof 
The proof is mainly based on the observation that the summands in both expressions in \eqref{eqn:Thetas} satisfy:
\begin{align*}
\prod_{\ell=0}^j (w+ m\ell) &= m^{j+1} \prod_{\ell=0}^j\left( \ell + \frac{w}{m}\right),
\end{align*}
and have the following natural lower and upper bounds, 
\begin{align}
\label{eqn:upper_and_Lower_bound}
m^{j+1} \prod_{\ell=0}^j\left( \ell + \left\lfloor\frac{w}{m}\right\rfloor \right)\le \prod_{\ell=0}^j (w+ m\ell) \le  m^{j+1} \prod_{\ell=0}^j\left( \ell + \left\lfloor\frac{w}{m}\right\rfloor +1\right).
\end{align}

Apply the upper bound in \eqref{eqn:upper_and_Lower_bound}, we have
\begin{align*}
\overline{\Theta}(k, w)=\sum_{j=k}^\infty \frac{(\g m)^j}{\prod_{\ell=0}^j(w+ m\ell) } & \ge \frac{1}{m}\sum_{j=k}^\infty\frac{\g ^j}{\prod_{\ell=0}^j\left( \ell + \lfloor\frac{w}{m}\rfloor +1 \right) } \\ &= \frac{\lfloor\frac{w}{m}\rfloor!}{m}\sum_{j=k}^\infty \frac{\g^j}{\left(j+\lfloor\frac{w}{m}\rfloor+1\right)!}\\ &= \frac{\lfloor\frac{w}{m}\rfloor!}{m\g^{\lfloor\frac{w}{m}\rfloor+1}e^{-\g}}\sum_{j=k}^\infty \frac{ e^{-\g} \g^{\left(j+\lfloor\frac{w}{m}\rfloor+1\right)}}{\left(j+\lfloor\frac{w}{m}\rfloor+1\right)!}\\ &=
\frac{1}{m\g \pi_{\lfloor\frac{w}{m}\rfloor}}
P\left[A_\g \ge k+\left\lfloor\frac{w}{m}\right\rfloor+1\right].
\end{align*}
Apply the lower bound in \eqref{eqn:upper_and_Lower_bound}, we have
\begin{align*}
\overline{\Theta}(k, w) = \sum_{j=k}^\infty \frac{(\g m)^j}{\prod_{\ell=0}^j(w+ m\ell) } & 
\le \frac{1}{m}\sum_{j=k}^\infty\frac{\g ^j}{\prod_{\ell=0}^j\left( \ell + \left\lfloor\frac{w}{m}\right\rfloor \right) } \\ 
&= \frac{\left\lfloor\frac{w}{m}-1\right\rfloor_+!}{m}\sum_{j=k}^\infty \frac{\g^j}{\left(j+\left\lfloor\frac{w}{m}\right\rfloor\right)!}
\quad \mbox{ (if $w<m$, $\lfloor\frac{w}{m}\rfloor=0$, and $\lfloor\frac{w}{m}-1\rfloor_+=0$)}
\\ 
    &= \frac{\left\lfloor\frac{w}{m}-1\right\rfloor_+!}{m\g^{\left\lfloor\frac{w}{m}\right\rfloor}e^{-\g} } \cdot \sum_{j=k}^\infty \frac{e^{-\g} \g^{\left(j+\left\lfloor\frac{w}{m}\right\rfloor\right)}}{\left(j+\left\lfloor\frac{w}{m}\right\rfloor\right)!}\\ 
&\le
\frac{1}{\max(1, \lfloor \frac{w}{m} \rfloor) \cdot 
m\pi_{\lfloor\frac{w}{m}\rfloor}}
\cdot P\left[A_\g \ge k+\left\lfloor\frac{w}{m}\right\rfloor\right],
\end{align*}
where the last inequality is because:
if $w<m$, then $\lfloor\frac{w}{m}\rfloor=0$, $\lfloor\frac{w}{m}-1\rfloor_+=0$, then
\begin{align}  \label{eqn:wm0}
\frac{\left\lfloor\frac{w}{m}-1\right\rfloor_+!}{m\g^{\left\lfloor\frac{w}{m}\right\rfloor}e^{-\g} } = \frac{1}{m \pi_0} =\frac{1}{m \pi_{\lfloor\frac{w}{m}\rfloor}} 
;
\end{align}
and if $w \ge m$, then $\lfloor \frac{w}{m} \rfloor \ge 1$, and
\begin{align}  \label{eqn:wm1}
\frac{\left\lfloor\frac{w}{m}-1\right\rfloor_+!}{m\g^{\left\lfloor\frac{w}{m}\right\rfloor}e^{-\g} } =\frac{1}{m \lfloor\frac{w}{m}\rfloor  \pi_{\lfloor\frac{w}{m}\rfloor}} 
\end{align}
Therefore, \eqref{eqn:barTheta} follows.

Similarly,
\begin{align*}
\Theta(k, w) = \sum_{j=0}^{k} \frac{(\g m)^j}{\prod_{\ell=0}^j(w+ m\ell) } & \ge \frac{1}{m}\sum_{j=0}^{k} \frac{\g ^j}{\prod_{\ell=0}^j\left( \ell + \lfloor\frac{w}{m}\rfloor+1 \right) } \\ &= \frac{\lfloor\frac{w}{m}\rfloor!}{m}\sum_{j=0}^{k}  \frac{\g^j}{\left(j+\lfloor\frac{w}{m}\rfloor+1\right)!}\\ &= \frac{\lfloor\frac{w}{m}\rfloor!}{m\g^{\lfloor\frac{w}{m}\rfloor+1}e^{-\g}}\sum_{j=0}^{k}  \frac{e^{-\g}\g^{\left(j+\lfloor\frac{w}{m}\rfloor+1\right)}}{\left(j+\lfloor\frac{w}{m}\rfloor+1\right)!}\\ &=
\frac{1}{m\g \pi_{\lfloor\frac{w}{m}\rfloor}}
P\left[\left\lfloor\frac{w}{m}\right\rfloor+1\le A_\g \le k+\left\lfloor\frac{w}{m}\right\rfloor+1\right].
\end{align*}

Apply the lower bound in \eqref{eqn:upper_and_Lower_bound}, we have
\begin{align*}
\Theta(k, w) = \sum_{j=0}^{k} \frac{(\g m)^j}{\prod_{\ell=0}^j(w+ m\ell) } & \le \frac{1}{m}\sum_{j=0}^{k} \frac{\g ^j}{\prod_{\ell=0}^j\left( \ell + \lfloor\frac{w}{m}\rfloor \right) } \\ 
&= \frac{\lfloor\frac{w}{m}-1\rfloor_+!}{m}\sum_{j=0}^{k}  \frac{\g^j}{\left(j+\lfloor\frac{w}{m}\rfloor\right)!}   \quad \mbox{ (if $w<m$, $\lfloor\frac{w}{m}\rfloor=0$ and $\lfloor\frac{w}{m}-1\rfloor_+=0$)}
\\ 
&= \frac{ \lfloor\frac{w}{m}-1\rfloor_+!}{m\g^{\lfloor\frac{w}{m}\rfloor}e^{-\g}}\sum_{j=0}^{k}  \frac{e^{-\g} \g^{\left(j+\lfloor\frac{w}{m}\rfloor\right)}}{\left(j+\lfloor\frac{w}{m}\rfloor\right)!}\\ 
&\le \frac{1}{\max(1, \lfloor \frac{w}{m} \rfloor) \cdot 
m \pi_{\lfloor\frac{w}{m}\rfloor}}
P\left[\left\lfloor\frac{w}{m}\right\rfloor\le A_\g \le k+\left\lfloor\frac{w}{m}\right\rfloor\right]
\end{align*}
where the last inequality is due to \eqref{eqn:wm0} and \eqref{eqn:wm1}.

Therefore, \eqref{eqn:Theta} follows.
\endproof

\smallskip
\section{Proof of Lemma~\ref{lem:smallW}} \label{sec:smallW}
\proof 
Since $w<\g m$, based on  \eqref{eqn:difference_rel_II} from Lemma~\ref{lem:difference}, we have
\begin{align} \label{eqn:negDf}
 - Df_h(w) = f_h(w) -f_h(w+m) 
 =  \left(1-\frac{w}{\g m}\right) f_h(w) + \frac{P({A}_\g\ge y)}{\g m} 
 \end{align}

Based on \eqref{eqn:SteinEq4}, we know
 \begin{align}  \label{eqn:LBUB}
- P({A}_\g < y) \overline{\Theta}(k_w, w)\le  f_h(w) \le P({A}_\g \ge y) \Theta(k_w-1, w),
 \end{align}
where $k_w=\lceil \frac{(my-w)^+}{m}\rceil = y - \lfloor\frac{w}{m}\rfloor$ since $y$ is an integer. 

Apply the upper bound of \eqref{eqn:LBUB} to \eqref{eqn:negDf}, we have
\begin{align}
-Df_h(w) &  \le  \left(1-\frac{w}{\g m}\right) 
P({A}_\g\ge y) \Theta(k_w-1, w)+ \frac{P({A}_\g\ge y)}{\g m} \nonumber \\
& \le P({A}_\g\ge y) \left( \left(1-\frac{w}{\g m}\right) 
\frac{
P\left(\left\lfloor\frac{w}{m}\right\rfloor\le A_\g \le k_w -1 +\left\lfloor\frac{w}{m}\right\rfloor\right)
}{ 
m\pi_{\lfloor\frac{w}{m}\rfloor}}
+ \frac{1}{\g m} 
\right) \nonumber \\
& \le {P({A}_\g\ge y)} \left( \left(1-\frac{w}{\g m}\right) 
\frac{
P\left(
A_\g < y\right)
}{
m\pi_{\lfloor\frac{w}{m}\rfloor}}
+ \frac{1}{\g m} 
\right)  \label{eqn:UB1}
\end{align}
where the second inequality is due to \eqref{eqn:Theta}, and the last inequality is because
$k_w -1 +\left\lfloor\frac{w}{m}\right\rfloor=y-1<y$.

Similarly, apply the lower bound of \eqref{eqn:LBUB} to \eqref{eqn:negDf}, we have
 \begin{align}
 - Df_h(w) &\ge 
-  \left(1-\frac{w}{\g m}\right) 
P({A}_\g < y) \hat\Theta(k_w, w)
+ \frac{P({A}_\g\ge y)}{\g m} \nonumber \\
&\ge 
-  \left(1-\frac{w}{\g m}\right) 
P({A}_\g < y) 
\frac{1}{
m\pi_{\lfloor\frac{w}{m}\rfloor}} 
P\left(A_\g \ge k_w+\left\lfloor\frac{w}{m}\right\rfloor\right)
+ \frac{P({A}_\g\ge y)}{\g m} \nonumber \\
& = -  \left(1-\frac{w}{\g m}\right) 
\frac{P({A}_\g < y) 
}{
m\pi_{\lfloor\frac{w}{m}\rfloor}}
P\left(A_\g \ge y\right)
+ \frac{P({A}_\g\ge y)}{\g m} 
\quad \quad \mbox{(since $k_w +\lfloor\frac{w}{m}\rfloor =y$)}
\nonumber\\
&\ge {P\left(A_\g \ge y\right)} \cdot 
\left[  -  \left(1-\frac{w}{\g m}\right) 
\frac{P({A}_\g < y) 
}{
m\pi_{\lfloor\frac{w}{m}\rfloor}}
- \frac{1}{m\g}
\right]   \label{eqn:LB1}
 \end{align}
where the second inequality is due to \eqref{eqn:barTheta}.

Combine \eqref{eqn:UB1} and \eqref{eqn:LB1}, we therefore have, for all $w< \g m$,
\begin{align}
\frac{| f_h(w) -f_h(w+m)|}{P({A}_\g\ge y)}
&\le  \Big| \left(1-\frac{w}{\g m}\right) 
\frac{P\left(A_\g < y\right)}{
m\pi_{\lfloor\frac{w}{m}\rfloor}}
+ 
\frac{1}{m\g} \Big| \nonumber\\
& \le
\frac{1}{m \pi_0 }+\frac{1}{m\g}. 
\label{eqn:w_smaller}
\end{align}
where the second inequality is due to the fact that 
the Poisson probabilities are increasing before the mean value $\g$ (since $\frac{w}{m} <\g$).

Therefore, \eqref{eqn:w_small} follows.
\endproof

\medskip
\section{Proof of Lemma \ref{lem:mediumW}}
\label{sec:BMminus}
Since $\g m \le w < my$,  $k_w=\lceil \frac{(my-w)}{m}\rceil
= y - \lfloor\frac{w}{m}\rfloor$ (since $y$ is an integer). 
 In this case, $k_w \ge 1$, and
 \begin{align} \label{eqn:k_wy}
 m\g \le m(y -1)  \le w+(k_w-1)m<my, 
\end{align}
 since $k_w$ is the minimum number of $m$-jumps, starting from $w$, to reach or go above $my$, and $y > \g +1$.
Based on Lemma~\ref{lem:grid_cocavity}, for any $j=1, ..., k_w-1$, 
\eqref{eqn:2nd_order_diff_k} can be rewritten as 
\begin{align*}
&  - Df_h(w+jm)+ (j+1)P({A}_\g \ge y)\cdot T_m(y)\\ 
&\ge  - Df_h(w+(j-1)m)+ jP({A}_\g \ge y)\cdot T_m(y)
\end{align*}
Thus,
\begin{align*} 
-Df_h(w+ (k_w-1)m) + k_w P({A}_\g \ge y) \cdot T_m(y) \ge
-Df_h(w) + P({A}_\g \ge y)\cdot T_m(y) > 0,
\end{align*}
where the second inequality is due to Lemma~\ref{lem:decreasing}.
Therefore, for any $m \g \le w < my$, we have,
\begin{align} 
\label{eqn:w_overall}
-Df_h(w)
\le -Df_h(w+ (k_w-1)m) + (k_w-1) P({A}_\g \ge y) \cdot T_m(y).
\end{align}
Let us calculate the first term of the RHS of \eqref{eqn:w_overall}. 
Apply (\ref{eqn:difference_rel_II}) we have, 
\begin{align}
\label{eqn:difference_m_grid_large}
- Df_h(w+ (k_w-1) m)  = &-f_h(w+ (k_w-1) m)  \left(\frac{w+ (k_w-1) m}{\g m}-1\right) + \frac{P({A}_\g\ge y)}{\g m}.
\end{align}
Note that $k_{w+(k_w-1)m} = 1$, hence, 
based on \eqref{eqn:split}, we have
\begin{align}
- f_h(w+ (k_w-1) m)= & 
-\frac{P({A}_\g\ge y)}{\mbox{\textcolor{black}{$w+ (k_w-1) m$}}}
+ P({A}_\g < y)\sum_{j=1}^\infty\frac{(\g m)^{j}}{\prod_{\ell=0}^{j} (w+ (k_w-1) m+m \ell) } \nonumber\\ 
\le & - \frac{P({A}_\g\ge y)}{\textcolor{black}{my}} + \sum_{j=1}^\infty\frac{(\g m)^{j}}{\prod_{\ell=0}^{j} m(\frac{w}{m}+ (k_w-1) + \ell) } \nonumber
\\ = & - \frac{P({A}_\g\ge y)}{\textcolor{black}{my}} + \frac{1}{m}\sum_{j=1}^\infty\frac{\g^{j}}{\prod_{\ell=0}^{j} ( \frac{w}{m}+(k_w-1) + \ell) } \nonumber
        \\ \le & - \frac{P({A}_\g\ge y)}{\textcolor{black}{\textcolor{black}{my}}} + \frac{1}{m}\sum_{j=1}^\infty\frac{\g^{j}}{\prod_{\ell=0}^{j} ( y-1 + \ell) }  \quad \mbox{(since $y =k_w + \lfloor\frac{w}{m}\rfloor$)}
        \nonumber \\ 
= &  -\frac{P({A}_\g\ge y)}{\textcolor{black}{my}} + \frac{( y- 2)! }{m \g^{ y- 1} e^{-\g}}\sum_{j=1}^\infty\frac{e^{-\g} \g^{ y + j-1}}{ ( y + j-1)! } \nonumber
\\ = & -\frac{P({A}_\g\ge y)}{\textcolor{black}{my}} + \frac{P(A_\g \ge y)}{m \g \pi_{ y- 2}}.\label{eqn:difference_m_grid_large_calc}
\end{align} 
Plug \eqref{eqn:difference_m_grid_large_calc} into \eqref{eqn:difference_m_grid_large}, 
we have, 
\begin{align}
&- Df_h(w+ (k_w-1) m)
\nonumber \\ 
\le &\left[-\frac{P({A}_\g\ge y)}{\textcolor{black}{my}} + \frac{P(A_\g \ge y)}{m\g \pi_{y- 2}}\right] \left(\frac{w+ (k_w-1) m}{\g m}-1\right) + \frac{P({A}_\g\ge y)}{\g m}\nonumber
\\ 
= & \left[ \frac{P({A}_\g\ge y)}{m\g \pi_{y- 2}}\right] \left(\frac{w+ (k_w-1) m}{\g m}-1\right) + \frac{P({A}_\g\ge y)}{\textcolor{black}{my}}
+ \frac{P({A}_\g\ge y)}{\g m} 
\left( - \frac{w+ (k_w-1) m}{\textcolor{black}{my}} +1\right)
\nonumber \\
\le &\left[ \frac{P(A_\g \ge  y)}{m\g \pi_{ y- 2}}\right] \left(\frac{y}{\g }-1\right) 
+ \frac{P({A}_\g\ge y)}{\textcolor{black}{my}}
+ \frac{P({A}_\g\ge y)}{\g m} 
\left( - \frac{m(y-1)}{\textcolor{black}{my}} +1\right)
\quad \quad \quad \mbox{(based on \eqref{eqn:k_wy})}
\nonumber \\
= & P(A_\g \ge  y) \left[ \frac{y-\g}{m \g^2 \pi_{ y- 2}} + \frac{1}{my}\left(1+ \frac{1}{\g}\right)\right] 
\label{eqn:w_big}
\end{align}

Combine \eqref{eqn:w_overall} and \eqref{eqn:w_big}, we then have, for any $m \g \le w < my$,
\begin{align*}
\frac{ f_h(w) -f_h(w+m)}{P({A}_\g\ge y)}=\frac{-Df_h(w)}{P({A}_\g\ge y)}\le & 
\frac{y-\g}{m\g^2\pi_{ y- 2}}  + \frac{1}{my}(1+ \frac{1}{\g}) 
+  (k_w-1)T_m(y) \\
\le & 
\frac{y-\g}{m\g^2\pi_{ y- 2}}  
+ \frac{1}{my}(1+ \frac{1}{\g}) +
\frac{y-\g}{m\g \lfloor\g\rfloor\pi_{ y- 2}}
\\
\le & 
\frac{2(y-\g)}{m\g\lfloor\g \rfloor \pi_{ y- 2}}  + \frac{1}{my}(1+ \frac{1}{\g}) \\
= & 
\frac{2(y-\g)}{m\lfloor\g \rfloor (y-1)\pi_{ y- 1}}  + \frac{1}{my}(1+ \frac{1}{\g}) \\
\le & 
\frac{2}{m\lfloor\g \rfloor\pi_{ y- 1}}  + \frac{1}{my}(1+ \frac{1}{\g}).
\end{align*}
where the second inequality is because $k_w-1 = y- \lfloor \frac{w}{m} \rfloor -1 \le y-\g$, since $\g \le \frac{w}{m}\le \lfloor \frac{w}{m}\rfloor +1$.

Therefore, \eqref{eqn:g1_bound_new} follows.
\endproof

\bigskip
\section{Proof of Lemma \ref{lem:Bdminus}} \label{sec:Bdminus}
From \eqref{eqn:SteinEq5}, we have for any $w>0$ and $d=1, ..., m-1$,
\begin{align*}
f_h(w) =&P({A}_\g \ge y) \overline{\Theta}(0,w) -  \overline{\Theta}(k_w, w), \\
f_h(w+d) =&P({A}_\g \ge y) \overline{\Theta}(0,w+d) -  \overline{\Theta}(k_{w+d}, w+d).
\end{align*}
Thus,
\begin{align} 
f_h(w) - f_h(w+d) 
&= P({A}_\g\ge y) \underbrace{\left( \overline{\Theta}(0,w) - \overline{\Theta}(0,w+d) \right)}_{\Delta_1} +
 \underbrace{\left(\overline{\Theta}(k_{w+d},w+d)-\overline{\Theta}(k_w,w)  \right)}_{\Delta_2}
\label{eqn:dfd}
\end{align}
It is easy to see that 
\begin{align} \label{eqn:F12}
|f_h(w) - f_h(w+d)| \le P({A}_\g\ge y) \Delta_1 + |\Delta_2|.
\end{align}
where $\Delta_1 :=\overline{\Theta}(0,w) - \overline{\Theta}(0,w+d)$ is always positive, and 
$\Delta_2:=\overline{\Theta}(k_{w+d},w+d)-\overline{\Theta}(k_w,w)$ depends on the values of $w$, $d$ and $k_{w+d}$.

Apply the bound \eqref{eqn:barTheta} to $\overline{\Theta}(0, w)$, we have
\begin{align}
0< \Delta_1 
\le  \overline{\Theta}(0,w) 
\le \frac{P(A_\g \ge \lfloor\frac{w}{m} \rfloor)
}{\max(1, \lfloor\frac{w}{m} \rfloor)
m \pi_{\lfloor\frac{w}{m}\rfloor}} 
\le \frac{1}{ \max(1, \lfloor\frac{w}{m} \rfloor) m \pi_{\lfloor\frac{w}{m}\rfloor}}.
\label{eqn:F1}
\end{align}

Meanwhile, 
\begin{align} \label{eqn:caseII}
|\Delta_2|\le \max(\overline{\Theta}(k_{w+d},w+d), \overline{\Theta}(k_w,w)).
\end{align}
Since $w < w+d < w+m$,  we have $k_w \ge k_{w+d} \ge k_{w+m}= k_w -1$.
Hence, either $k_{w+d}=k_w$, or $k_{w+d}=k_w-1$.
\begin{itemize}
\item
$k_{w+d}=k_w$: $\overline{\Theta}(k_{w+d},w+d) 
= \overline{\Theta}(k_{w},w+d) 
<\overline{\Theta}(k_w,w)$, hence, $|\Delta_2| < \overline{\Theta}(k_w,w)$. Based on \eqref{eqn:barTheta},
\begin{align}
\overline{\Theta}(k_w,w) 
\le \frac{P\left[A_\g \ge k_w + \left\lfloor\frac{w}{m}\right\rfloor\right] 
}{\max(1, \lfloor\frac{w}{m} \rfloor)
m\pi_{\lfloor\frac{w}{m}\rfloor}}
= \frac{P(A_\g \ge y)}{\max(1, \lfloor\frac{w}{m} \rfloor) m \pi_{\lfloor\frac{w}{m}\rfloor}}.
\label{eqn:F2}
\end{align} 
\item
$k_{w+d}=k_w-1$:  In this case,  $k_{w+d}= y - \lfloor\frac{w+d}{m}\rfloor =k_w -1= y - \lfloor\frac{w}{m}\rfloor -1
$ implies that $\lfloor\frac{w+d}{m}\rfloor=\lfloor\frac{w}{m}\rfloor+1$,
and $k_w -1+\lfloor\frac{w+d}{m}\rfloor  =y$.
Apply the bound \eqref{eqn:barTheta} to $\overline{\Theta}(k_{w+d}, w+d)$, we have
\begin{align} 
 \overline{\Theta}(k_{w+d},w+d) 
\le  \frac{P\left[A_\g \ge k_w -1 + \left\lfloor\frac{w+d}{m}\right\rfloor\right] 
}{\max(1, \lfloor\frac{w+d}{m} \rfloor)
m\pi_{\lfloor\frac{w+d}{m}\rfloor}} 
=  \frac{P(A_\g \ge y)}{(\lfloor\frac{w}{m} \rfloor+1)
m\pi_{\lfloor\frac{w}{m}\rfloor+1}} 
=  \frac{P(A_\g \ge y)}{
m\g \pi_{\lfloor\frac{w}{m}\rfloor}}.  \label{eqn:caseIIb}
\end{align}

Combine \eqref{eqn:caseII}-\eqref{eqn:caseIIb}, we then have:
\begin{align} \label{eqn:caseIIc}
|\Delta_2|  \le {P(A_\g \ge y)} \cdot 
\max\left(\frac{1}{m\g \pi_{\lfloor\frac{w}{m}\rfloor}}, 
\frac{1}{ \max(1, \lfloor\frac{w}{m} \rfloor) m \pi_{\lfloor\frac{w}{m}\rfloor}}\right)
\end{align}
\end{itemize}

Combine \eqref{eqn:F12}, \eqref{eqn:F1}, and \eqref{eqn:caseIIc}, we then have 
\begin{align}
\frac{|f_h(w+d) -f_h(w)| }{P({A}_\g\ge y)} 
\le \frac{2}{m \min\{ \pi_0,  \lfloor \g \rfloor \pi_{y-1}\}},
\end{align}
because if $w < m\g$, then $0 \le \lfloor\frac{w}{m}\rfloor < \g$, 
$\pi_{\lfloor\frac{w}{m}\rfloor}\ge \pi_0$, and we have
$\frac{1}{m\g \pi_{\lfloor\frac{w}{m}\rfloor}}\le \frac{1}{m\pi_{0}}$, and 
$\frac{1}{\max(1, \lfloor\frac{w}{m} \rfloor) m \pi_{\lfloor\frac{w}{m}\rfloor}}
\le \frac{1}{m\pi_{0}}$;
if $m \g \le w < my$, then 
$\g \le \lfloor \frac{w}{m} \rfloor \le y-1$, 
$\pi_{\lfloor\frac{w}{m}\rfloor}\ge \pi_{y-1}$, and we have
$\frac{1}{m\g \pi_{\lfloor\frac{w}{m}\rfloor}}\le \frac{1}{m
\lfloor \g \rfloor \pi_{y-1}}$, and 
\[\frac{1}{\max(1, \lfloor\frac{w}{m} \rfloor) m \pi_{\lfloor\frac{w}{m}\rfloor}}
=\frac{1}{m \lfloor\frac{w}{m} \rfloor \pi_{\lfloor\frac{w}{m}\rfloor}} 
=\frac{1}{m \g \pi_{\lfloor\frac{w}{m}\rfloor-1}}
\le \frac{1}{m \lfloor \g \rfloor \pi_{y-2}}
\le \frac{1}{m \lfloor \g \rfloor \pi_{y-1}}
\]
where the first inequality is based on \eqref{eqn:Piy_2}.
\endproof

\section{Proof of Proposition~\ref{pro:LX_conjecture}} \label{sec:conj}
\proof 
From~\eqref{eqn:C1-} and \eqref{eqn:C1+}, we can see that, the comparison between ${\mathfrak C}_{1-}(\g,k)$ and ${\mathfrak C}_{1+}(\g,k)$ can be reduce to the compaison of $\frac{F(k-2)}{F(k-1)} \cdot \frac{\g}{k-1}$ and $\frac{{\bar F}(k+1)}{{\bar F}(k)} \cdot \frac{k}{\g}$. 
Examine their ratio, we have, 
\begin{align*}
\frac{\frac{F(k-2)}{F(k-1)} \cdot \frac{\g}{k-1}}{\frac{{\bar F}(k+1)}{{\bar F}(k)} \cdot \frac{k}{\g}}
 =& \frac{F(k-2){\bar F}(k)}{F(k-1){\bar F}(k+1)}\cdot \frac{\g^2}{k(k-1)} 
\\ =& \frac{F(k-2)({\bar F}(k+1)+\pi_k)}{(F(k-2)+\pi_{k-1}){\bar F}(k+1)}\cdot \frac{\g^2}{k(k-1)} 
\\ =& \frac{F(k-2){\bar F}(k+1)+F(k-2)\pi_k}{F(k-2){\bar F}(k+1)+\pi_{k-1}{\bar F}(k+1)}\cdot \frac{\g^2}{k(k-1)} 
\\ =& \frac{1+ \frac{\pi_k}{{\bar F}(k+1)}}{1+ \frac{\pi_{k-1}}{F(k-2)}}\cdot \frac{\g^2}{k(k-1)}.
\end{align*}
Meanwhile, we know,
\begin{align*}
\frac{\pi_k}{{\bar F}(k+1)}\le \frac{\pi_k}{\pi_{k+1}}=\frac{k+1}{\g}.
\end{align*}
Therefore, 
\begin{align*}
\frac{1+ \frac{\pi_k}{{\bar F}(k+1)}}{1+ \frac{\pi_{k-1}}{F(k-2)}}\cdot \frac{\g^2}{k(k-1)}\le \left(1+ \frac{k+1}{\g}\right)\cdot \frac{\g^2}{k(k-1)}.
\end{align*}
It is easy to verify that the desired inequality holds if the condition is satisfied. 
\endproof

\end{appendix}

\end{document}